\documentclass[10pt]{article}

\def \N {{\bf {N}}}
\def \R {{\bf {R}}}

\def \M {{\cal M}}
\def \m {{\bf m}}

\title{\bf On spectral multiplicities of  Gaussian actions}
\author{V.V.Ryzhikov}

\date{30.06. 2014}
\begin{document}
\LARGE
\maketitle
\begin{abstract}{\large For any set $M\subset \N$ there are 
 mixing Gaussian  automorphisms and 
   non-mixing  Gaussian  automorphisms (as well as automorphisms which are disjoint from any Gaussian action)
with continuous  spectrum  for which    $ M\cup\{\infty\}$ is the set of their spectral multiplicities. We show also that for a  Gaussian
flow  $\{G_t\}$  the  spectral multiplicities of $G_t$, \ $t>0$, could be different.  }
\end{abstract}

\section{ Spectral multiplicities of  some Gaussian automorphisms} \rm Let  $Exp(U)=\bigoplus_{n=1}^{\infty} U^{\odot n}$, where  $U$ is  \bf a unitary operator with simple continuous  singular \rm spectrum.
\vspace{4mm}

\bf Lemma.  \it  There are $z_m$, $m\in \N$, $|z_m|=1$, such that the  operators 
$z_m U$ are pairwise spectrally disjoint. \rm

{\large (Proof. $U$ has singular spectrum = its spectral  measure $\sigma$ is singular.  Let  $\sigma_z$, $|z|=1$, is defined  by the condition
$$\sigma_z(A)=\sigma (zA).$$
The measure $\sigma$ is singular iff for almost all $z$ with respect to Lebesgue measure on the unit circle  the measures
$\sigma_z$ and $\sigma$ are disjoint.) } 

Denote  $\m U = \bigoplus_{n=1}^{m} U$.
Given $M\subset \N$, we set
 $V=  \bigoplus_{m\in M}\m (z_m U),$ where 
all  $z_m$ are from Lemma.
\vspace{4mm}

\bf Assertion 1. \it     If 
a) the product $U^{\odot 2}$ has absolutely continuous
spectrum, 
or 

b) the operators  $U^{\odot n}$ have homogeneous spectra of infinite multiplicity as  $n>1$, and  all operators $U^{\odot n}$
are pairwise disjoint,

then  the set $\M(Exp(V))$ of spectral multiplicities
of the  operator $Exp(V)$  is
$$ M\cup\{\infty\}.$$
\rm

Proof.  It is not hard to see that $$Exp(V)=  V  \bigoplus  \bigoplus_{i=1}^\infty W_{(i)}, $$
where in a) the operator $W$ has Lebesgue spectrum,
and in b) $W$ has singular spectrum.
\vspace{4mm}

\bf Realization of a). \rm The examples of singular measures  $\sigma$  (considered  as a spectral measure for a unitary operator)
with the absolutely continuous convolution $\sigma\ast\sigma$ are classical (Schaeffer-Salem measures, 1939-1943, see \cite{P}).

\bf Realization of b). \rm  There is an operator $U$ of simple spectrum  such that

  $2U^{n_i}\to_w I$ (providing the  operators $U^{\odot n}$
to be pairwise disjoint), and

  $U^3$ is spectrally isomorphic to $U\oplus U\oplus U$  (this implies that 
 the operators  $U^{\odot n}$ have homogeneous spectra of infinite multiplicity as $n>1$).

 We apply here  the operator
induced by an exponential "self-similar" 2-adic Chacon map, see  arXiv:1311.4524.

\vspace{4mm}
There is well-known correspondence between  a unitary operator $U$ with continuous spectrum
and  an ergodic  Gaussian automorphism (of a probability space) $G(U)$  which is spectrally isomorphic to the operator
$$\bigoplus_{n=0}^{\infty} U^{\odot n},$$
where $U^{\odot 0}=1$ denotes one-dimension identity
operator (acting on the space of constants).
\vspace{4mm}

Thus, from the above facts we get the following result,
adding a string to the list of realizable sets (the sets of spectral multiplicities for ergodic automorphisms).
\vspace{4mm}

\bf Theorem 1. \it    For any set $M\subset \N$ there are 
 mixing Gaussian  automorphisms and 
   non-mixing    automorphisms
with singular spectrum  for which    $ M\cup\{\infty\}$ is the set of their spectral multiplicities. \rm
\vspace{1mm}
\\
\bf Spectral multiplicities of  non-Gaussian automorphisms.

\vspace{1mm}
\bf Theorem 1.1 \it    For any set $M\subset \N$ there is a weakly mixing  automorphism $T$ which is disjoint with any 
  Gaussian  automorphism and $\M(T)= M\cup\{\infty\}$. \rm
\vspace{1mm}

For this we consider
$T=T_1\times T_2\times \dots$, where $T_i$  have homogeneous 
 spectra with multiplicities from $M$ and all pairwise convolutions of their 
spectral measures are Lebesgue.

We recall that all sets in the form $ M\cup\{1\}$ and $ M\cup\{2\}$ are realizable  via   special skew products (Ageev, Danilenko, Katok,
Kwiatkowski, Lemanczyk et al,  see \cite{D}).  
\vspace{4mm}

\section{Spectral multiplicities of   automorphisms in the Gaussian  flow} \rm  
Some results on non-usual multiplicities of powers for a weakly mixing automorphism see in \cite{R}.\rm  Any Gaussian automorphism is embedded in a Gaussian flow.
We prove that
 \it  there are  Gaussian
flows  $\{G_t\}$  with non-constant sets $\M (G_t) $ (of the  spectral multiplicities for automorphisms  $G_t$, \ $t>0$).  \rm

We consider an operator $U=V\oplus (-V)$ with  simple spectrum, where $V$ is disjoint with  the powers  $V^{\odot n}$, $n\geq 2$, 
  and the corresponding Gaussian flow $G_t$, $G_1=G(U)$.
Then we get \footnote{\large This answers the question 11 from \cite{D}}
\vspace{4mm}

\bf Theorem 2. \it There is $U$ such that  $\M(G(U)) =\{1,\dots\}$ and  $\M(G(U^2)) =\{2,\dots\}\subset 2\N$. \rm
\vspace{4mm}

Now we change an operator setting $U$ as in section 1, the case  b. The operator $U$ has simple spectrum (it is induced by a rank-one transformation). Its power $U^{3^k}$
is isomorphic to the direct sum of $3^k$ copies of $U$. All symmetric powers $U^{\odot 2}, \dots,$  have infinite spectral multiplicity. Thus, we obtain
\vspace{2mm}

\bf Theorem 3. \it There is $U$ such that   $\M(G(U^{3^k})) =\{3^k,\infty\}.$ \rm 
\vspace{1mm}

However  a general simple fact shows that $\M(G_t)$ is constant  almost everywhere.
\vspace{4mm}

\bf Assertion 2. \it  Let $\{U_t\}$ be a unitary flow with simple singular spectrum.  Then for almost all $t$ the operators $U_t$ have 
simple spectra as well.  Generally, if the spectral measure of a flow is singular, then    $\M(U_t)=
\M(\{U_t\})$ for almost all $t$.  \rm
\vspace{4mm}

By  use of combinations of operators
in the form $ \bigoplus_{m=1}^{n}z_m U$ we get a variety 
of flows  with non-constant multiplicity function $\M(G_t)$.

\vspace{4mm}
\bf Theorem 4. \it 1. For any prime $p$ and any $m<p$ there is a Gaussian
flow $\{G_t\}$ such that   $\M(G_1)=\{1,\infty\} $ and $\M(G_p)=\{m,\infty\} $.

2.  There is a Gaussian
flow $\{G_t\}$ such that  for any finite integer set $M$   for some
$t=t(M)>0$  $\M(G_t)=\{1,\infty\}\cup M $. 

3. 
For any   set $M\subset \N$ there is a Gaussian
flow $\{G_t\}$ with singular spectrum such that  $\M(\{G_t\})=\{1,\infty\}$, $\M(G_1)=M\cup \{\infty\}$, and   any interval $I\subset \R$ contains  an infinite set  $\{t_k\}$  for which multiplicitiy sets $\M(G_{t_k})$ are pairwise different.

4. Let $P$ denote the set of primes. For any function
$m:P\to \N$, $m(p)\leq p$, there is a  Gaussian
flow $\{G_t\}$ such that 
$$ \M(G_n)=\{\infty\}\cup \{m(p_1), m(p_2),\dots, m(p_k)\},$$
where $n=p_1^{d_1}p_2^{d_2}\dots p_k^{d_k}$,  the  
factorization of $n$ into a product of  primes. 
\vspace{4mm}

\rm 

We thank J.-P. Thouvenot  for useful discussions.
\large

 \end{document}